\documentclass{amsart}
\usepackage{amsthm}
\usepackage{amsmath}
\usepackage{epsfig}
\usepackage[usenames]{color} % pickup color
\usepackage{ulem} % pickup sout

\newtheorem{remark}{Remark}
\newtheorem*{theorem*}{Theorem}%%%%%%%%%%%use with \usepackage{amsthm}
\begin{document}

\numberwithin{equation}{section}

\title{Locally CR Spherical Three Manifolds}
\author{Howard Jacobowitz}
\address{Department of Mathematical Sciences,
Rutgers University, Camden, NJ 08102}
\email{jacobowi@camden.rutgers.edu}
\date{March 21, 2013}
\begin{abstract}%%%%%%%%%%%%%%%%%%%%   should precede maketitle
Every open and orientable three manifold has a CR structure which is locally equivalent to the standard CR structure on $S^3$.
\end{abstract}
\keywords{CR, locally spherical}
%\subjclass{32V05} gave 1991 subject class!
\thanks{2010 Mathematics Subject Classification: 32V05}
\maketitle
\def\spc{strictly pseudo-convex}
\def\psh{plurisubharmonic}
\def\spsh{strictly plurisubharmonic}
\def\nbd{neighborhood}
\def\iff{if and only if }
\def\comtanM{C\otimes T(M)\ }
\def\T10{T^{1,0}C^N}
\def\CN{C^N\ }
\def\comtanO{C\otimes T(\Omega)}
\def\dpartial{\bar{\partial}}

The standard CR structure on $S^3$ is the one induced by the usual complex structure on $R^4$.  
  A CR structure on a three dimensional manifold  is said to be locally spherical if for each point there exist an open neighborhood $U$ and a CR diffeomorphism of $U$ onto an open subset of $S^3$.  There are many results about compact locally spherical three manifolds  (for instance, \cite{BE},  \cite{CL}, \cite{FG}, \cite{I}, \cite{L}) which show that these manifolds are special.  The situation is quite different for open manifolds.  This easy observation seems not to be in the literature.  (There are some interesting constructions on special open manifolds, particularly in connection with complex hyperbolic geometry \cite{F}, \cite{S}.)
\begin{theorem*}
Every open and orientable three dimensional manifold admits a locally spherical CR structure.
\end{theorem*}
\begin{remark}
The CR structure may be taken to be $C^\omega$.
\end{remark}
If we remove one point from $S^3$ and choose any diffeomorphism of this new manifold to $R^3$, we end up with a CR structure on $R^3$.  (Usually this diffeomorphism is taken so that we end up with the CR structure induced on the hyperquadric $\{ (z,w)\in C^2 \ |\ \Im{w}=\lvert{z^2}\rvert\}$.)  Any map
\[
F:M\to R^3
\]
which is a local diffeomorphism may be used to pull this CR structure on $R^3  $ back to $M$.  The result is a locally spherical CR structure on $M$.  The existence of some $F$ is a well-known topological fact, first proved by Whitehead \cite{W}.  See also Phillips \cite{P} , Corollary 8.2 and the references cited there.  
\begin{remark} The conclusion of the theorem also holds for any parallelizable manifold of odd dimension.
\end{remark}

\end{document}